\documentclass [12pt]{article}
\usepackage{amssymb}
\def\lanbox{\hbox{$\, \vrule height 0.25cm width 0.25cm depth 0.01cm \,$}}

\usepackage{color}

\begin{document}

\centerline{\Large\bf Solvability of
some integro-differential equations}

\centerline{\Large\bf with the bi-Laplacian and transport}

\bigskip

\centerline{Vitali Vougalter}

\bigskip

\centerline{Department of Mathematics, University of Toronto,
Toronto, Ontario, M5S 2E4, Canada}

\centerline{e-mail: vitali@math.toronto.edu}

\bigskip
\bigskip
\bigskip

\noindent {\bf Abstract.}
We demonstrate the existence in the sense of sequences of solutions for
some integro-differential type problems involving the drift term
and the square of the Laplace operator, on the
whole real line or on a finite interval with periodic boundary conditions
in the corresponding $H^{4}$ spaces. Our argument is based on the fixed point
technique when the elliptic equations contain fourth order differential
operators with and without the Fredholm property. It is established that,
under the reasonable technical conditions, the convergence in $L^{1}$ of the
integral kernels yields the existence and convergence in $H^{4}$ of the
solutions.

\bigskip
\bigskip



\bigskip
\bigskip
\bigskip
\bigskip

\section{Introduction}

 We recall that a linear operator $L$, which acts from a Banach
 space $E$ into another Banach space $F$ satisfies the Fredholm
 property if its image is closed, the dimension of its kernel and
 the codimension of its image are finite. Consequently, the
 equation $Lu=f$ is solvable if and only if $\phi_i(f)=0$ for a
 finite number of functionals $\phi_i$ from the dual space $F^*$.
 Such properties of the Fredholm operators are broadly used in many
 methods of the linear and nonlinear analysis.

 Elliptic problems in bounded domains with a sufficiently smooth
 boundary possess the Fredholm property if the ellipticity
 condition, proper ellipticity and Lopatinskii conditions are
 fulfilled (see e.g. \cite{Ag}, \cite{E09}, \cite{LM}, \cite{Volevich}).
 This is the main result of the theory of linear elliptic equations.
 In the situation of unbounded domains, such conditions may
 not be sufficient and the Fredholm property may not be satisfied.
 For instance, the Laplace operator, $Lu = \Delta u$, in $\mathbb R^d$
 fails to satisfy the Fredholm property when considered in
 H\"older spaces, $L : C^{2+\alpha}(\mathbb R^d) \to C^{\alpha}(\mathbb
 R^d)$, or in Sobolev spaces,  $L : H^2(\mathbb R^d) \to L^2(\mathbb
 R^d)$.

 Linear elliptic problems in unbounded domains satisfy the
 Fredholm property if and only if, in addition to the conditions
 stated above, the limiting operators are invertible (see \cite{V11}).
 In some simple cases, the limiting operators can be constructed
 explicitly. For example, if
 $$
 L u = a(x) u'' + b(x) u' + c(x) u , \;\;\; x \in \mathbb R ,
 $$
 where the coefficients of the operator have the limits at infinities,

 $$ a_\pm =\lim_{x \to \pm \infty} a(x) , \;\;\;
 b_\pm =\lim_{x \to \pm \infty} b(x) , \;\;\;
 c_\pm =\lim_{x \to \pm \infty} c(x) , $$
 the limiting operators are given by:

 $$ L_{\pm}u = a_\pm u'' + b_\pm u' + c_\pm u . $$
 Since the coefficients here are constants, the essential spectrum of
 the operator, that is the set of complex numbers $\lambda$ for
 which the operator $L-\lambda$ does not satisfy the Fredholm
 property, can be explicitly found by means of the Fourier
 transform:
 $$
 \lambda_{\pm}(\xi) = -a_\pm \xi^2 + b_\pm i\xi + c_\pm , \;\;\;
 \xi \in \mathbb R .
 $$
The invertibility of limiting operators is equivalent to the condition that 
 the essential spectrum does not contain the origin.

 In the situation of general elliptic equations, the analogous assertions hold true.
 The Fredholm property is satisfied if the the origin does not belong to the essential spectrum 
 or if the limiting operators are invertible. However,
 these conditions may not be written explicitly.

 In the case of the non-Fredholm operators the usual solvability
 relations may not be applicable and the solvability conditions
 are, in general, not known. There are certain classes of operators
 for which the solvability relations are obtained. We illustrate
 them with the following example. Consider the equation
 \begin{equation}
 \label{int1}
 Lu \equiv \Delta u + a u = f
 \end{equation}
 in $\mathbb R^d, \ d\in {\mathbb N}$, where $a>0$ is a constant.
 The operator $L$ coincides with its limiting operators. The
 homogeneous problem admits a nontrivial bounded solution. Thus, the
 Fredholm property is not satisfied. However, because our operator
 has constant coefficients, we can use the Fourier transform and
 find the solution explicitly. The solvability conditions can be
 formulated as follows. If $f \in L^2(\mathbb R^d)$ and 
 $xf \in L^1(\mathbb R^d)$, then there
 exists a unique solution of this problem in $H^2(\mathbb R^d)$ if and
 only if
 $$  
 \Bigg(f(x),\frac{e^{ipx}}{(2\pi)^{\frac{d}{2}}}\Bigg)_{L^2(\mathbb R^d)}=0, \quad 
 p\in S_{\sqrt{a}}^{d} \quad a.e.
 $$
 (see \cite{VV103}). Here $S_{\sqrt{a}}^{d}$ denotes the sphere
 in $\mathbb R^d$ of radius $\sqrt{a}$ centered at the origin.
 Therefore, though the operator fails to satisfy the Fredholm property,
 the solvability relations are formulated similarly. However,
 this similarity is only formal since the range of our operator is
 not closed.

 In the situation of the operator involving a scalar potential,
 $$
 L u \equiv \Delta u + a(x) u = f ,
 $$
 the Fourier transform is not directly applicable. Nevertheless, the solvability
 conditions in ${\mathbb R}^{3}$ can be obtained by a rather sophisticated 
 application of the theory of self-adjoint
 operators (see \cite{VV08}). As before, the solvability relations are
 formulated in terms of the orthogonality to the solutions of the homogeneous
 adjoint
 problem. There are several other examples of linear elliptic
 non-Fredholm operators for which the solvability
 conditions can be derived (see ~\cite{EV21}, ~\cite{V11}, ~\cite{VKMP02}, ~\cite{VV11}
 ~\cite{VV08}, ~\cite{VV10}, ~\cite{VV103}).

 Solvability relations play a crucial role in the analysis of
 nonlinear elliptic equations. In the case of the operators without the Fredholm
 property, in spite of a certain progress in the understanding of the linear
 problems, there exist only few examples where nonlinear non-Fredholm
 operators are analyzed (see ~\cite{DMV08},
 ~\cite{VV103}, ~\cite{VV21}, ~\cite{VV211}). The work ~\cite{E10} deals with
 the studies of the finite and infinite dimensional attractors for evolution
 equations of mathematical physics. The large time behavior of solutions
 of a class of fourth-order parabolic equations defined on unbounded domains
 using the Kolmogorov $\varepsilon$-entropy as a measure was covered
 in ~\cite{EP07}. In ~\cite{EZ01} the authors discuss the attractor for a
 nonlinear reaction-diffusion system in an unbounded domain in
 ${\mathbb R}^{3}$. The articles ~\cite{GS05} and ~\cite{RS01} are dedicated to
 the understanding of the Fredholm and properness properties of quasilinear
 elliptic systems of second order and of the operators of this kind on
 ${\mathbb R}^{N}$. In ~\cite{GS10} the authors demonstrate the exponential
 decay and study the Fredholm properties in second order quasilinear elliptic
 systems. The work ~\cite{E18} is a systematic study of a dynamical
 systems approach to investigating the symmetrization and stabilization
 properties
 of nonnegative solutions of nonlinear elliptic equations in asymptotically
 symmetric unbounded domains.
 In the current article we deal with another class of stationary
 nonlinear equations, for which the Fredholm property may not be satisfied:
\begin{equation}
\label{id1}
-\frac{d^{4}u}{dx^{4}}+ b\frac{du}{dx}+au + 
\int_{\Omega}G(x-y)F(u(y),y)dy = 0, \quad x\in \Omega.
\end{equation}
Here $\Omega$ is a domain on ${\mathbb R}, \ a\geq 0, \ b\in {\mathbb R}, \ b\neq 0$ are the constants. 
For the simplicity of the presentation we restrict ourselves to the one dimensional situation
(the multidimensional case is more technical). The solvability of
the integro-differential problem with the transport and influx/efflux terms
on the real line was considered in ~\cite{VV211}. The method used there
worked for the powers of the fractional negative Laplacian
$\displaystyle{0<s<\frac{1}{4}}$.
In the population dynamics the
integro-differential equations describe the models with the intra-specific
competition and the nonlocal consumption of resources (see e.g.
\cite{ABVV10}, \cite{BNPR09}, ~\cite{VV13}).  Let us use the explicit form of
the solvability conditions and establish the existence of solutions of our
nonlinear problem. The studies of the solutions of the integro-differential
equations with the drift term are relevant to the understanding of the
emergence and propagation of patterns in the theory of speciation
(see ~\cite{VV13}).
The solvability of the linear problem containing the Laplace operator with
the transport term was discussed in ~\cite{VV10}, see also ~\cite{BHN05}.
In the situation when the drift term is trivial, namely when $b=0$, the equation
analogous to (\ref{id1}) was treated in ~\cite{VV21} with the constant $a>0$.
In article ~\cite{AMP14} the authors consider the wave systems with an infinite
number of localized traveling waves. Standing lattice solitons in the discrete NLS equation with saturation
were covered in ~\cite{AKLP19}. Diffusion-driven 
blow-up for a nonlocal Fisher-KPP type model was discussed in ~\cite{KLS23}.
Markov process representations for polyharmonic functions were presented in ~\cite{K10}.

%%%%%%%%%%%%%%%%%%%%%%%%%%%%%%%%%%%%%%%%%%%%%%%%%%%%%%%%%%%%%%%%

\setcounter{equation}{0}

\section{Formulation of the results}

The nonlinear part of equation (\ref{id1}) will satisfy the
following regularity assumptions.

\bigskip

\noindent
{\bf Assumption 1.} {\it Function $F(u,x): {\mathbb R}\times
\Omega \to {\mathbb R}$ is satisfying the Caratheodory condition 
(see ~\cite{K64}), such that
\begin{equation}
\label{ub1}
|F(u,x)|\leq k|u|+h(x) \quad for \quad u\in {\mathbb R}, \ x\in \Omega
\end{equation}
with a constant $k>0$ and $h(x):\Omega\to {\mathbb R}^{+}, \quad
h(x)\in L^{2}(\Omega)$. Moreover, it is a Lipschitz continuous function,
such that
\begin{equation}
\label{lk1}
|F(u_{1},x)-F(u_{2},x)|\leq l |u_{1}-u_{2}| \quad for \quad any \quad
 u_{1,2}\in{\mathbb R}, \quad x\in \Omega
\end{equation}
with a constant $l>0$.}

\bigskip

The work ~\cite{BO86} deals with the solvability of a local elliptic
problem in a bounded domain in ${\mathbb R}^{N}$. The nonlinear function contained there
was allowed to have a sublinear growth.
In order to establish the solvabity of equation (\ref{id1}), we will use
the auxiliary problem
\begin{equation}
\label{ae1}
\frac{d^{4}u}{dx^{4}}-b\frac{du}{dx}-au=\int_{\Omega}G(x-y)F(v(y),y)dy, \quad x\in \Omega,
\end{equation}
where $a\geq 0, \ b\in {\mathbb R}, \ b\neq 0$ are the constants.
We designate
\begin{equation}
\label{ip}  
(f_{1}(x),f_{2}(x))_{L^{2}(\Omega)}:=\int_{\Omega}f_{1}(x)\bar{f_{2}}(x)dx,
\end{equation}
with a slight abuse of notations when the functions involved in  (\ref{ip}) do not belong to $L^{2}(\Omega)$,
like for example those contained in orthogonality relation
(\ref{or1}) below. Clearly, if $f_{1}(x)\in L^{1}(\Omega)$ and
$f_{2}(x)$ is bounded, then the integral in the right side of definition (\ref{ip})
makes sense.
In the first part of the work we consider the case on
the whole real line,  $\Omega={\mathbb R}$, so that the appropriate Sobolev
space is equipped with the norm
\begin{equation}
\label{h2n}  
\|u\|_{H^{4}({\mathbb R})}^{2}:=\|u\|_{L^{2}({\mathbb R})}^{2}+
\Bigg\|\frac{d^{4}u}{dx^{4}}\Bigg\|_{L^{2}({\mathbb R})}^{2}.
\end{equation}
The main issue for the problem above is that in the absence
of the transport term we deal with the self-adjoint, non-Fredholm  operator
$$
\frac{d^{4}}{dx^{4}}-a: H^{4}({\mathbb R})\to L^{2}({\mathbb R}), \
a\geq 0,
$$
which is the obstacle to solve our equation (see ~\cite{VV21}). Note that in article ~\cite{VV21}
it was assumed that the constant $a>0$ since the problem becomes more singular when $a$ vanishes.
The similar cases but
in linear problems, both self- adjoint and non self-adjoint
involving the differential operators without the Fredholm property
have been studied broadly in recent years (see
~\cite{V11}, ~\cite{VKMP02}, ~\cite{VV08}, ~\cite{VV10}, ~\cite{VV103}).
 However, the situation is different when
the constant in the transport term $b$ is nontrivial. The operator
\begin{equation}
\label{lab}
L_{a, \ b}:=\frac{d^{4}}{dx^{4}}-b\frac{d}{dx}-a: \quad
H^{4}({\mathbb R})\to L^{2}({\mathbb R}),
\end{equation}
where the constants $a\geq 0$ and $b\in {\mathbb R}, \ b\neq 0$ contained in the left side
of problem (\ref{ae1}) is non-selfadjoint. By means of the standard Fourier
transform, it can be easily derived that the essential spectrum of
this operator $L_{a, \ b}$ is given by
$$
\lambda_{a, \ b}(p)=p^{4}-a-ibp, \quad p\in {\mathbb R},
$$
(see ~\cite{V11}).
Obviously, when the constant $a>0$ the operator $L_{a, \ b}$ satisfies
the Fredholm property, since the origin does not belong to its
essential spectrum. But when $a=0$, our
operator $L_{a, \ b}$ fails to satisfy the Fredholm property because  its essential spectrum
 contains the origin. We establish that under the reasonable
technical conditions problem (\ref{ae1}) defines a map $t_{a, b}:
H^{4}({\mathbb R})\to H^{4}({\mathbb R})$ with the constants
$a\geq 0, \ b\in {\mathbb R}, \ b\neq 0$, which is a strict contraction.

\bigskip

\noindent
{\bf Theorem 1.}  {\it Let $\Omega={\mathbb R}, \ G(x): {\mathbb R}\to
{\mathbb R}, \ G(x)\in L^{1}({\mathbb R})$ and Assumption 1  is valid.

\medskip

\noindent
I) If  $a>0, \ b\in {\mathbb R}, \ b\neq 0$, we assume that
$2\sqrt{\pi}N_{a, \ b}l<1$ with $N_{a, \ b}$ defined in (\ref{Na}).
Then the map $v \mapsto t_{a, b}v=u$ on $H^{4}({\mathbb R})$ defined by equation
(\ref{ae1}) has a unique fixed point $v_{a, b}$. This is the only
solution of problem (\ref{id1}) in $H^{4}({\mathbb R})$.

\medskip

\noindent
II) Let $a=0, \ b\in {\mathbb R}, \ b\neq 0$. Assume that
$xG(x)\in L^{1}({\mathbb R})$,
orthogonality relation (\ref{or1}) holds
and $2\sqrt{\pi}N_{0, \ b}l<1$.
Then the map $t_{0, b}v=u$ on $H^{4}({\mathbb R})$ defined by equation
(\ref{ae1}) possesses  a unique fixed point $v_{0, b}$, which is
the only solution of problem (\ref{id1}) in
$H^{4}({\mathbb R})$.

\medskip

\noindent
In both cases I and II the fixed point
$v_{a, b}, \ a\geq 0, \ \ b\in {\mathbb R}, \ b\neq 0$ does not vanish identically in ${\mathbb R}$
provided the intersection of supports of the Fourier transforms
of functions $supp\widehat{F(0,x)}\cap supp \widehat{G}$ is a set of nonzero
Lebesgue measure on the real line.}

\bigskip

Note that in the case of this theorem  when $a>0$, as distinct
from Theorem 1 of ~\cite{VV21} describing the equation without
the transport term, the orthogonality conditions are not needed.
Let us introduce the sequence of
approximate equations related to problem (\ref{id1}) on ${\mathbb R}$,
namely
\begin{equation}
\label{id1m}
-\frac{d^{4}u_{m}}{dx^{4}}+ b\frac{du_{m}}{dx}+au_{m} + 
\int_{-\infty}^{\infty}G_{m}(x-y)F(u_{m}(y),y)dy = 0,
\end{equation}
where $a\geq 0, \ b\in {\mathbb R}, \ b\neq 0$ are the constants and
$m\in {\mathbb N}$.
The sequence of kernels $\displaystyle{\{G_{m}(x)\}_{m=1}^{\infty}}$ tends
to $G(x)$ as $m\to \infty$ in the corresponding function spaces specified below.
We establish that, under the appropriate technical conditions, each of problems
(\ref{id1m}) admits a unique solution
$u_{m}(x)\in H^{4}({\mathbb R})$, limiting equation (\ref{id1})
has a unique solution $u(x)\in H^{4}({\mathbb R})$, and
$u_{m}(x)\to u(x)$ in $H^{4}({\mathbb R})$ as $m\to \infty$. This is the
so-called {\it existence of solutions in the sense of sequences}. In this
situation, the solvability conditions can be formulated for the iterated kernels
$G_{m}$. They yield the convergence of the kernels in terms of the Fourier
transforms (see the Appendix) and consequently the convergence of the
solutions (Theorems 2, 4).
The analogous ideas in the sense of the standard Schr\"odinger type operators
were used in ~\cite{VV131}. Our second main result is as
follows.

\bigskip

\noindent
{\bf Theorem 2.}  {\it Let $\Omega={\mathbb R}, \ m\in {\mathbb N}, \
G_{m}(x): {\mathbb R}\to {\mathbb R}, \ G_{m}(x)\in L^{1}({\mathbb R})$ are
such that $G_{m}(x)\to G(x)$ in $L^{1}({\mathbb R})$ as $m\to \infty$ and
Assumption 1  holds.

\medskip

\noindent
I) When $a>0, \ b\in {\mathbb R}, \ b\neq 0$, assume that
$$  
2\sqrt{\pi}N_{a, \ b, \ m}l\leq 1-\varepsilon
$$
for all $m\in {\mathbb N}$ with some fixed $0<\varepsilon<1$ and
$N_{a, \ b, \ m}$ introduced in (\ref{Nam}).
Then each equation (\ref{id1m}) possesses a unique solution
$u_{m}(x)\in H^{4}({\mathbb R})$, and limiting problem (\ref{id1})
has a unique solution $u(x)\in H^{4}({\mathbb R})$.

\medskip

\noindent
II) For $a=0, \ b\in {\mathbb R}, \ b\neq 0$, assume that
$xG_{m}(x)\in L^{1}({\mathbb R}), \ xG_{m}(x)\to xG(x)$ in $L^{1}({\mathbb R})$
as $m\to \infty$, orthogonality relation (\ref{or1m}) is valid and 
$$ 
2\sqrt{\pi}N_{0, \ b, \ m}l\leq 1-\varepsilon
$$
for all $m\in {\mathbb N}$ with a certain fixed $0<\varepsilon<1$.
Then each problem (\ref{id1m}) admits a unique solution
$u_{m}(x)\in H^{4}({\mathbb R})$, and limiting equation (\ref{id1})
possesses a unique solution $u(x)\in H^{4}({\mathbb R})$.

\medskip

\noindent
In both cases I and II, we have $u_{m}(x)\to u(x)$ in $H^{4}({\mathbb R})$
as $m\to \infty$.

\medskip

\noindent
The unique solution $u_{m}(x)$ of each equation (\ref{id1m}) is nontrivial
provided that the intersection of supports of the Fourier transforms
of functions $supp\widehat{F(0,x)}\cap supp \widehat{G}_{m}$ is a set of nonzero
Lebesgue measure on the real line. Analogously, the unique solution $u(x)$
of limiting problem (\ref{id1}) does not vanish identically in  ${\mathbb R}$ if
$supp\widehat{F(0,x)}\cap supp \widehat{G}$ is a set of nonzero Lebesgue
measure on the real line.}

\bigskip

The second part of the work is devoted to the studies of the analogous equation on the
finite interval $\Omega=I:=[0, \ 2\pi]$ with periodic boundary conditions.
The corresponding function space is given by
$$
H^{4}(I)=\{u(x):I\to {\mathbb R} \ | \ u(x), u''''(x)\in L^{2}(I), \quad
u(0)=u(2\pi), \quad u'(0)=u'(2\pi), 
$$
$$
u''(0)=u''(2\pi), \quad u'''(0)=u'''(2\pi) \}.
$$
For the technical purposes, we introduce the following auxiliary constrained
subspace
\begin{equation}
\label{H0}
H_{0}^{4}(I)=\{u(x)\in H^{4}(I) \ | \ (u(x),1)_{L^{2}(I)}=0 \},
\end{equation}
which is a Hilbert space as well (see e.g. Chapter 2.1 of ~\cite{HS96}).
Our goal is to demonstrate that problem (\ref{ae1}) in such case defines a map
$\tau_{a, b}$ on the above mentioned
spaces with the constants $a\geq 0, \ b\in {\mathbb R}, \ b\neq 0$. This map
will be a strict contraction under the given technical
conditions.

\bigskip

\noindent
{\bf Theorem 3.} {\it Let
$\Omega=I, \ G(x): I \to {\mathbb R}, \ G(x)
\in C(I), \ G(0)=G(2\pi), \ F(u,0)=F(u, 2\pi)$
for $u\in {\mathbb R}$ and Assumption 1 holds.

\medskip

\noindent
I) Let $a>0, \ b\in {\mathbb R}, \ b\neq 0$. Assume that
$2\sqrt{\pi}{\cal N}_{a, \ b}l<1$, where
${\cal N}_{a, \ b}$ is introduced in (\ref{Nab}).
Then the map $v \mapsto \tau_{a, b}v=u$ on
$H^{4}(I)$ defined by equation (\ref{ae1}) has a unique fixed point
$v_{a, b}$, the only solution of problem (\ref{id1}) in $H^{4}(I)$.

\medskip

\noindent
II) If $a=0, \ b\in {\mathbb R}, \ b\neq 0$, we assume that orthogonality
relation (\ref{oc1}) is valid and
$2\sqrt{\pi}{\cal N}_{0, \ b}l<1$. Then the map $\tau_{0, b}v=u$ on $H^{4}_{0}(I)$
defined by problem (\ref{ae1}) possesses a unique fixed point $v_{0, b}$, the
only solution of equation (\ref{id1}) in $H^{4}_{0}(I)$.

\medskip

\noindent
In both cases I and II the fixed point
$v_{a, b}, \ a\geq 0, \ b\in {\mathbb R}, \ b\neq 0$ is nontrivial
on the interval $I$
provided the Fourier coefficients $G_{n}F(0, x)_{n}\neq 0$ for some
$n\in {\mathbb Z}$.}

\bigskip

\noindent
{\bf Remark 1.} {\it We use the constrained subspace $H^{4}_{0}(I)$ 
in case II) of the theorem, such that the Fredholm operator
$\displaystyle{\frac{d^{4}}{dx^{4}}-b\frac{d}{dx}: H^{4}_{0}(I)\to L^{2}(I)}$
has the trivial kernel.}

\bigskip

Let us establish the solvability in the sense of sequences for our
integro-differential problem on the interval $I$. Consider the sequence of approximate equations
similarly to the situation on the whole real line with $m\in {\mathbb N}$,
namely
\begin{equation}
\label{id1mi}
-\frac{d^{4}u_{m}}{dx^{4}}+ b\frac{du_{m}}{dx}+au_{m} + 
\int_{0}^{2\pi}G_{m}(x-y)F(u_{m}(y),y)dy = 0.
\end{equation}
Here $a\geq 0, \ b\in {\mathbb R}, \ b\neq 0$ are the constants. Our final
main statement is as follows.

\bigskip

\noindent
{\bf Theorem 4.} {\it Let
$\displaystyle{\Omega=I, \ m\in {\mathbb N}, \ G_{m}(x): I\to{\mathbb R}, \
G_{m}(x)\in C(I)}$,
such that
$$  
G_{m}(x)\to G(x) \quad in \quad C(I) \quad
as \quad m\to \infty, 
$$
$G_{m}(0)=G_{m}(2\pi), \ F(u, 0)=F(u, 2\pi)$
for $u\in {\mathbb R}$. Let Assumption 1 be valid.

\medskip

\noindent
I) For $a>0, \ b\in {\mathbb R}, \ b\neq 0$, assume that 
$$  
2\sqrt{\pi}{\cal N}_{a, \ b, \ m}l\leq 1-\varepsilon
$$
for all $m\in {\mathbb N}$ with some fixed $0<\varepsilon<1$ and
${\cal N}_{a, \ b, \ m}$ introduced in (\ref{Nabmc}).
Then each
problem (\ref{id1mi}) possesses a unique solution $u_{m}(x)\in H^{4}(I)$ and 
limiting equation (\ref{id1}) admits a unique solution $u(x)\in H^{4}(I)$.

\medskip

\noindent
II) When $a=0, \ b\in {\mathbb R}, \ b\neq 0$, assume that the orthogonality 
relation (\ref{or1mi}) holds and
$$  
2\sqrt{\pi}{\cal N}_{0, \ b, \ m}l\leq 1-\varepsilon
$$
for all $m\in {\mathbb N}$ with a certain fixed $0<\varepsilon<1$. Then each
equation (\ref{id1mi}) admits a unique solution $u_{m}(x)\in H^{4}_{0}(I)$ and 
limiting problem (\ref{id1}) has a unique solution $u(x)\in H^{4}_{0}(I)$.

\medskip

\noindent
In both cases I and II we have $u_{m}(x)\to u(x)$ as $m\to \infty$ in the norms
in $H^{4}(I)$ and $H^{4}_{0}(I)$ respectively.

\medskip

\noindent
The unique solution $u_{m}(x)$ of each equation (\ref{id1mi}) is nontrivial
on the interval $I$
provided that the Fourier coefficients $G_{m, n}F(0, x)_{n}\neq 0$ for a certain
$n\in {\mathbb Z}$. Analogously, the unique solution $u(x)$ of
limiting problem (\ref{id1}) does not vanish identically on the interval $I$ if
$G_{n}F(0, x)_{n}\neq 0$ for some $n\in {\mathbb Z}$.}

\bigskip

\noindent
{\bf Remark 2.} {\it In the manuscript we work with the real
valued functions by means of the conditions imposed on $F(u,x), \ G_{m}(x)$ and $G(x)$
contained in the integral terms of the iterated and limiting equations discussed
above.}

\bigskip

\noindent
{\bf Remark 3.} {\it The significance of Theorems 2 and 4 above is the continuous
dependence of the solutions with respect to the integral kernels.}      

%%%%%%%%%%%%%%%%%%%%%%%%%%%%%%%%%%%%%%%%%%%%%%%

\setcounter{equation}{0}

\section{The Whole Real Line Case}

\bigskip

{\it Proof of Theorem 1.} First we suppose that in the case of
$\Omega={\mathbb R}$ for some $v\in H^{4}({\mathbb R})$ there
exist two solutions $u_{1,2}\in H^{4}({\mathbb R})$ of equation
(\ref{ae1}). Then their difference
$w(x):=u_{1}(x)-u_{2}(x)\in  H^{4}({\mathbb R})$ will be a solution of the homogeneous
problem
$$
\frac{d^{4}w}{dx^{4}}-b\frac{dw}{dx}-aw=0.
$$
Since the operator
$L_{a, \ b}: H^{4}({\mathbb R})\to L^{2}({\mathbb R})$ defined in (\ref{lab})
does not have any nontrivial zero modes,
$w(x)$ vanishes identically on  the real line.

Let us choose an arbitrary $v(x)\in H^{4}({\mathbb R})$ and apply the standard
Fourier transform (\ref{ft}) to both sides of (\ref{ae1}). This gives us
\begin{equation}
\label{f1}
\widehat{u}(p)=\sqrt{2\pi}\frac{\widehat{G}(p)\widehat{f}(p)}{p^{4}-a-ibp}, \quad
p^{4}\widehat{u}(p)=\sqrt{2\pi}\frac{p^{4}\widehat{G}(p)\widehat{f}(p)}
{p^{4}-a-ibp}.
\end{equation}
Here $\widehat{f}(p)$ stands for the Fourier image of $F(v(x),x)$. Evidently,
the upper bounds
$$
|\widehat{u}(p)|\leq \sqrt{2\pi}N_{a, \ b}|\widehat{f}(p)|  \quad and
\quad |p^{4}\widehat{u}(p)|\leq \sqrt{2\pi}N_{a, \ b}|\widehat{f}(p)|
$$
hold.
We have $N_{a, \ b}<\infty$ by virtue of Lemma A1 of the Appendix without
any orthogonality conditions if $a>0$ and under orthogonality relation
(\ref{or1}) for $a=0$. This allows us to derive the estimate from above
on the norm
$$
\|u\|_{H^{4}({\mathbb R})}^{2}=
$$
\begin{equation}
\label{h4nub}
\|\widehat{u}(p)\|_{L^{2}({\mathbb R})}^{2}+
\|p^{4}\widehat{u}(p)\|_{L^{2}({\mathbb R})}^{2}\leq 4\pi N_{a, \ b}^{2}\int_{-\infty}^{\infty}|\widehat{f}(p)|^{2}dp=
4\pi N_{a, \ b}^{2}\|F(v(x),x)\|_{L^{2}({\mathbb R})}^{2}.
\end{equation}
Let us recall  inequality (\ref{ub1}) of Assumption 1 above. Hence, the right side of  (\ref{h4nub}) is finite for
$v(x)\in L^{2}({\mathbb R})$. This means that for an arbitrary
$v(x)\in H^{4}({\mathbb R})$ there exists a unique solution
$u(x)\in H^{4}({\mathbb R})$ of equation (\ref{ae1}). Its Fourier image
is given by (\ref{f1}) and the map
$t_{a, b}:  H^{4}({\mathbb R})\to  H^{4}({\mathbb R})$ is well defined.

This enables us to choose arbitrarily $v_{1,2}(x)\in H^{4}({\mathbb R})$,
so that their images $u_{1,2}:=t_{a, b}v_{1,2}\in H^{4}({\mathbb R})$.
By means of (\ref{ae1}), we have
\begin{equation}
\label{aer1}
\frac{d^{4}u_{1}}{dx^{4}}-b\frac{du_{1}}{dx}-au_{1}=\int_{-\infty}^{\infty}
G(x-y)F(v_{1}(y),y)dy,
\end{equation}
\begin{equation}
\label{aer2}
\frac{d^{4}u_{2}}{dx^{4}}-b\frac{du_{2}}{dx}-au_{2}=\int_{-\infty}^{\infty}
G(x-y)F(v_{2}(y),y)dy.
\end{equation}
We apply the standard Fourier transform (\ref{ft}) to both sides of
problems (\ref{aer1}) and (\ref{aer2}). This yields
$$
\widehat{u}_{1}(p)=\sqrt{2\pi}\frac{\widehat{G}(p)\widehat{f}_{1}(p)}
{p^{4}-a-ibp}, \quad
p^{4}\widehat{u}_{1}(p)=\sqrt{2\pi}\frac{p^{4}\widehat{G}(p)\widehat{f}_{1}(p)}
{p^{4}-a-ibp},
$$
$$
\widehat{u}_{2}(p)=\sqrt{2\pi}\frac{\widehat{G}(p)\widehat{f}_{2}(p)}
{p^{4}-a-ibp}, \quad
p^{4}\widehat{u}_{2}(p)=\sqrt{2\pi}\frac{p^{4}\widehat{G}(p)\widehat{f}_{2}(p)}
{p^{4}-a-ibp}.
$$
Here $\widehat{f}_{1}(p)$ and  $\widehat{f}_{2}(p)$ denote the Fourier transforms
of $F(v_{1}(x),x)$  and $F(v_{2}(x),x)$ respectively.
Obviously, the upper bounds
$$
|\widehat{u}_{1}(p)-\widehat{u}_{2}(p)|\leq \sqrt{2\pi}N_{a, \ b}
|\widehat{f}_{1}(p)-\widehat{f}_{2}(p)|,  \quad
|p^{4}\widehat{u}_{1}(p)-p^{4}\widehat{u}_{2}(p)|\leq \sqrt{2\pi}
N_{a, \ b}|\widehat{f}_{1}(p)-\widehat{f}_{2}(p)|
$$
are valid.
This enables us to derive the inequality for the norms
$$
\|u_{1}-u_{2}\|_{H^{4}({\mathbb R})}^{2}\leq 4\pi N_{a, \ b}^{2} \int_{-\infty}^{\infty}|\widehat{f}_{1}(p)-\widehat{f}_{2}(p)|^{2}dp=
4\pi N_{a, \ b}^{2}
\|F(v_{1}(x),x)-F(v_{2}(x),x)\|_{L^{2}({\mathbb R})}^{2}.
$$
Note that $v_{1,2}(x)\in H^{4}({\mathbb R})\subset
L^{\infty}({\mathbb R})$ due to the Sobolev embedding.
Let us recall condition (\ref{lk1}). This gives us
$$
\|F(v_{1}(x),x)-F(v_{2}(x),x)\|_{L^{2}({\mathbb R})}\leq l \|v_{1}(x)-v_{2}(x)\|_{L^{2}({\mathbb R})},
$$
such that
\begin{equation}
\label{Tab}  
\|t_{a, b}v_{1}-t_{a, b}v_{2}\|_{H^{4}({\mathbb R})}\leq 2
\sqrt{\pi}N_{a, \ b}l\|v_{1}-v_{2}\|_{H^{4}({\mathbb R})}.
\end{equation}
The constant in the right side of estimate (\ref{Tab}) is less than one as we assume.
According to the Fixed Point Theorem,
there exists a unique function
$v_{a, b}\in H^{4}({\mathbb R})$ having the property $t_{a, b}v_{a, b}=v_{a, b}$.
This is the only solution of equation (\ref{id1}) in
$H^{4}({\mathbb R})$. Suppose $v_{a, b}(x)$ vanishes identically in ${\mathbb R}$.
This will contradict to the given condition that the Fourier images
of $G(x)$ and $F(0,x)$ are nontrivial on a set of nonzero Lebesgue
measure on the real line.
\hfill\lanbox

\bigskip

We turn our attention to establishing the existence in the sense of
sequences of the solutions for our integro-differential problem on ${\mathbb R}$.

\bigskip

\noindent
{\it Proof of Theorem 2.} According to the result of Theorem 1 above, each
problem (\ref{id1m}) admits a unique solution
$u_{m}(x)\in H^{4}({\mathbb R}), \ m\in {\mathbb N}$. Limiting equation
(\ref{id1}) has a unique solution $u(x)\in H^{4}({\mathbb R})$ by virtue
of Lemma A2 below along with  Theorem 1. We apply the standard
Fourier transform (\ref{ft}) to both sides of (\ref{id1}) and (\ref{id1m})
and arrive at
\begin{equation}
\label{ump}
\widehat{u}(p)=\sqrt{2\pi}\frac{\widehat{G}(p)\widehat{\varphi}(p)}{p^{4}-a-ibp},
\quad \widehat{u}_{m}(p)=\sqrt{2\pi}\frac{\widehat{G}_{m}(p)
\widehat{\varphi}_{m}(p)}{p^{4}-a-ibp}, \quad m\in {\mathbb N}.
\end{equation}  
Here $\widehat{\varphi}(p)$ and $\widehat{\varphi}_{m}(p)$ designate the
Fourier images of $F(u(x), x)$ and $F(u_{m}(x), x)$ respectively. Clearly,
$$
|\widehat{u}_{m}(p)-\widehat{u}(p)|\leq \sqrt{2\pi}\Bigg\|\frac{\widehat{G}_{m}
(p)}{p^{4}-a-ibp}-\frac{\widehat{G}(p)}{p^{4}-a-ibp}\Bigg\|_{L^{\infty}
({\mathbb R})}|\widehat{\varphi}(p)|+
$$
$$
\sqrt{2\pi}\Bigg\|\frac{\widehat{G}_{m}(p)}{p^{4}-a-ibp}\Bigg\|_{L^{\infty}
({\mathbb R})}|\widehat{\varphi}_{m}(p)-\widehat{\varphi}(p)|.
$$
Therefore,
$$
\|u_{m}-u\|_{L^{2}({\mathbb R})}\leq \sqrt{2\pi}\Bigg\|\frac{\widehat{G}_{m}
(p)}{p^{4}-a-ibp}-\frac{\widehat{G}(p)}{p^{4}-a-ibp}\Bigg\|_{L^{\infty}
({\mathbb R})}\|F(u(x), x)\|_{L^{2}({\mathbb R})}+
$$
$$
\sqrt{2\pi}\Bigg\|\frac{\widehat{G}_{m}(p)}{p^{4}-a-ibp}\Bigg\|_{L^{\infty}
({\mathbb R})}\|F(u_{m}(x), x)-F(u(x), x)\|_{L^{2}({\mathbb R})}.
$$
We recall condition (\ref{lk1}) of Assumption 1. Hence,
\begin{equation}
\label{Fumu}  
\|F(u_{m}(x), x)-F(u(x), x)\|_{L^{2}({\mathbb R})}\leq l\|u_{m}(x)-u(x)\|_{L^{2}
({\mathbb R})}.
\end{equation}
Obviously,
$u_{m}(x), \ u(x)\in H^{4}({\mathbb R})\subset L^{\infty}({\mathbb R})$ via
the Sobolev embedding. Thus, we obtain
$$
\|u_{m}(x)-u(x)\|_{L^{2}({\mathbb R})}\Bigg\{1-\sqrt{2\pi}\Bigg\|\frac{\widehat{G}_{m}(p)}{p^{4}-a-ibp}\Bigg\|_{L^{\infty}({\mathbb R})}l\Bigg\}\leq
$$
$$
\sqrt{2\pi}
\Bigg\|\frac{\widehat{G}_{m}(p)}{p^{4}-a-ibp}-\frac{\widehat{G}(p)}
{p^{4}-a-ibp}\Bigg\|_{L^{\infty}({\mathbb R})}\|F(u(x), x)\|_{L^{2}({\mathbb R})}.
$$
Let us recall inequality (\ref{2rnabml}) for $a>0$ and (\ref{2rnabml0}) when $a=0$. This yields
$$
\|u_{m}(x)-u(x)\|_{L^{2}({\mathbb R})}\leq \frac{\sqrt{2\pi}}{\varepsilon}
\Bigg\|\frac{\widehat{G}_{m}(p)}{p^{4}-a-ibp}-\frac{\widehat{G}(p)}
{p^{4}-a-ibp}\Bigg\|_{L^{\infty}({\mathbb R})}\|F(u(x), x)\|_{L^{2}({\mathbb R})}.
$$
By virtue of condition (\ref{ub1}) of Assumption 1, we have
$F(u(x), x)\in L^{2}({\mathbb R})$ for $u(x)\in H^{4}({\mathbb R})$. This means that
\begin{equation}
\label{umuc}
u_{m}(x)\to u(x), \quad m\to \infty   
\end{equation}
in $L^{2}({\mathbb R})$ due to the result of Lemma A2 of the Appendix.
Note that
$$
p^{4}\widehat{u}(p)=\sqrt{2\pi}\frac{p^{4}\widehat{G}(p)\widehat{\varphi}(p)}
{p^{4}-a-ibp},\quad p^{4}\widehat{u}_{m}(p)=\sqrt{2\pi}\frac{p^{4}\widehat{G}_{m}(p)
\widehat{\varphi}_{m}(p)}{p^{4}-a-ibp}, \quad m\in {\mathbb N},
$$
such that
$$
|p^{4}\widehat{u}_{m}(p)-p^{4}\widehat{u}(p)|\leq \sqrt{2\pi}
\Bigg\|\frac{p^{4}\widehat{G}_{m}(p)}{p^{4}-a-ibp}-\frac{p^{4}\widehat{G}(p)}
{p^{4}-a-ibp}\Bigg\|_{L^{\infty}({\mathbb R})}|\widehat{\varphi}(p)|+
$$
$$
\sqrt{2\pi}\Bigg\|\frac{p^{4}\widehat{G}_{m}(p)}{p^{4}-a-ibp}\Bigg\|_{L^{\infty}
({\mathbb R})}|\widehat{\varphi}_{m}(p)-\widehat{\varphi}(p)|.
$$
Let us use (\ref{Fumu}) to derive 
$$
\Bigg\|\frac{d^{4}u_{m}}{dx^{4}}-\frac{d^{4}u}{dx^{4}}\Bigg\|_{L^{2}({\mathbb R})}\leq
\sqrt{2\pi}
\Bigg\|\frac{p^{4}\widehat{G}_{m}(p)}{p^{4}-a-ibp}-\frac{p^{4}\widehat{G}(p)}
{p^{4}-a-ibp}\Bigg\|_{L^{\infty}({\mathbb R})}\|F(u(x), x)\|_{L^{2}({\mathbb R})}+
$$
$$
\sqrt{2\pi}\Bigg\|\frac{p^{4}\widehat{G}_{m}(p)}{p^{4}-a-ibp}\Bigg\|_{L^{\infty}
({\mathbb R})}l\|u_{m}(x)-u(x)\|_{L^{2}({\mathbb R})}.
$$
We recall the result of Lemma A2 of the Appendix along with (\ref{umuc}). Therefore,
$\displaystyle{\frac{d^{4}u_{m}}{dx^{4}}\to \frac{d^{4}u}{dx^{4}}}$ in
$L^{2}({\mathbb R})$ as $m\to \infty$. According to definition (\ref{h2n}) of the norm,
we arrive at $u_{m}(x)\to u(x)$ in $H^{4}({\mathbb R})$ as $m\to \infty$.

Suppose that the unique solution of problem (\ref{id1m})
considered above  $u_{m}(x)\equiv 0$ for some $m\in {\mathbb N}$. This will give a contradiction
to our condition that the Fourier transforms of $G_{m}(x)$ and $F(0, x)$ are nontrivial
on a set of nonzero Lebesgue measure on the real line. The analogous reasoning
holds for the unique solution $u(x)$ of limiting equation (\ref{id1}).
\hfill\lanbox

 %%%%%%%%%%%%%%%%%%%%%%%%%%%%%%%%%%%%%%%%%%%%%%%

 \setcounter{equation}{0}

 \section{The Problem on the Finite Interval}

\bigskip

\noindent
{\it Proof of Theorem 3.} We present the proof of the theorem in the case
of $a>0$. When the constant   
$a=0$, the reasoning will be analogous. If $a=0$, we will need to use the
constrained subspace (\ref{H0}) instead of $H^{4}(I)$. The non-selfadjoint
operator contained in the left side of equation (\ref{ae1}) is given by
\begin{equation}
\label{clab}
{\cal L}_{a, \ b}:=\frac{d^{4}}{dx^{4}}-b\frac{d}{dx}-a: \quad H^{4}(I)\to
L^{2}(I).     
\end{equation}
Its set of eigenvalues is 
\begin{equation}
\label{spclab}  
{\cal \lambda}_{a, \ b}(n)=n^{4}-a-ibn, \quad n\in {\mathbb Z}
\end{equation}
and its square integrable eigenfunctions are the standard Fourier harmonics
$\displaystyle{\frac{e^{inx}}{\sqrt{2\pi}}, \ n\in {\mathbb Z}}$. Clearly,
the eigenvalues of the the operator ${\cal L}_{a, \ b}$ are simple, as distinct
from the similar case without the transport term, when the eigenvalues
corresponding to $n\neq 0$ are two-fold degenerate (see ~\cite{VV21}).
Thus, ${\cal L}_{a, \ b}$ satisfies the Fredholm property.

Let us first suppose that for a certain $v(x)\in H^{4}(I)$ there exist two
solutions $u_{1,2}(x)\in H^{4}(I)$ of equation (\ref{ae1}) with $\Omega=I$.
Then the function $w(x):=u_{1}(x)-u_{2}(x)\in  H^{4}(I)$ will solve the
homogeneous problem
$$
\frac{d^{4}w}{dx^{4}}-b\frac{dw}{dx}-aw=0.
$$
Evidently, the operator ${\cal L}_{a, \ b}: \ H^{4}(I)\to L^{2}(I)$ introduced above
does not have any nontrivial zero modes. Therefore,  $w(x)$ vanishes identically in $I$.

We choose an arbitrary $v(x)\in H^{4}(I)$, apply the Fourier
transform (\ref{fti}) to equation (\ref{ae1}) considered on the interval $I$
and arrive at
\begin{equation}
\label{f2}
u_{n}=\sqrt{2\pi}\frac{G_{n}f_{n}}{n^{4}-a-ibn}, \quad n^{4}u_{n}=\sqrt{2\pi}
\frac{n^{4}G_{n}f_{n}}{n^{4}-a-ibn}, \quad n\in {\mathbb Z},
\end{equation}
where $f_{n}:=F(v(x),x)_{n}$. This allows us to obtain the estimates from above
$$
|u_{n}|\leq \sqrt{2\pi}{\cal N}_{a, \ b}|f_{n}|, \quad |n^{4}u_{n}|\leq \sqrt{2\pi}
{\cal N}_{a, \ b}|f_{n}|.
$$
Obviously, ${\cal N}_{a, \ b}<\infty$ under the given conditions according to
Lemma A3 of the Appendix. Hence, we derive
$$
\|u\|_{H^{4}(I)}^{2}=
$$
\begin{equation}
\label{uh4iub}
\sum_{n=-\infty}^{\infty}|u_{n}|^{2}+
\sum_{n=-\infty}^{\infty}|n^{4}u_{n}|^{2}\leq 4\pi {\cal N}_{a, \ b}^{2}\sum_{n=-\infty}^{\infty}|f_{n}|^{2}
= 4\pi {\cal N}_{a, \ b}^{2}\|F(v(x),x)\|_{L^{2}(I)}^{2}.
\end{equation}
Recall condition (\ref{ub1}) of Assumption 1. Thus, the right side of (\ref{uh4iub}) is finite for $v(x)\in L^{2}(I)$.
Therefore, for an arbitrarily chosen $v(x)\in H^{4}(I)$ there exists a unique
$u(x)\in H^{4}(I)$, which satisfies equation (\ref{ae1}). Its
Fourier transform is given by (\ref{f2}), so the map
$\tau_{a, b}: H^{4}(I)\to H^{4}(I)$ in the case I) of the theorem is well
defined.

Consider arbitrary $v_{1,2}(x)\in H^{4}(I)$. Their images under the map
mentioned above $u_{1,2}:=\tau_{a, b}v_{1,2}\in H^{4}(I)$. According to (\ref{ae1}),
we have
\begin{equation}
\label{ae2pi1}
\frac{d^{4}u_{1}}{dx^{4}}-b\frac{du_{1}}{dx}-au_{1}=\int_{0}^{2\pi}G(x-y)
F(v_{1}(y),y)dy,
\end{equation}
\begin{equation}
\label{ae2pi2}
\frac{d^{4}u_{2}}{dx^{4}}-b\frac{du_{2}}{dx}-au_{2}=\int_{0}^{2\pi}G(x-y)
F(v_{2}(y),y)dy.
\end{equation}
Let us apply Fourier transform (\ref{fti}) to both sides of
(\ref{ae2pi1}) and (\ref{ae2pi2}). This gives us
$$
u_{1, n}=\sqrt{2\pi}\frac{G_{n}f_{1, n}}{n^{4}-a-ibn}, \quad u_{2, n}=\sqrt{2\pi}
\frac{G_{n}f_{2, n}}{n^{4}-a-ibn},
$$
$$
n^{4}u_{1, n}=\sqrt{2\pi}\frac{n^{4}G_{n}f_{1, n}}{n^{4}-a-ibn}, \quad n^{4}u_{2, n}=
\sqrt{2\pi}\frac{n^{4}G_{n}f_{2, n}}{n^{4}-a-ibn}, \quad n\in {\mathbb Z}.
$$
Here $f_{j, n}:=F(v_{j}(x),x)_{n},\ j=1,2$. Thus,
$$
|u_{1, n}-u_{2, n}|\leq \sqrt{2\pi}{\cal N}_{a, \ b}|f_{1, n}-f_{2, n}|, \quad
|n^{4}(u_{1, n}-u_{2, n})|\leq \sqrt{2\pi}{\cal N}_{a, \ b}|f_{1, n}-f_{2, n}|,
$$
such that
$$
\|u_{1}-u_{2}\|_{H^{4}(I)}^{2}= \sum_{n=-\infty}^{\infty}
|u_{1, n}-u_{2, n}|^{2}+
\sum_{n=-\infty}^{\infty}|n^{4}(u_{1, n}-u_{2, n})|\leq
$$
$$
 4\pi {\cal N}_{a, \ b}^{2}\sum_{n=-\infty}^{\infty}|f_{1, n}-f_{2, n}|^{2} =4\pi {\cal N}_{a, \ b}^{2}\|F(v_{1}(x),x)-F(v_{2}(x),x)\|_{L^{2}(I)}^{2}.
$$
Clearly, $v_{1,2}(x)\in H^{4}(I)\subset L^{\infty}(I)$ via the Sobolev
embedding. By means of condition (\ref{lk1}), we easily derive
$$
\|F(v_{1}(x),x)-F(v_{2}(x),x)\|_{L^{2}(I)}\leq l\|v_{1}(x)-v_{2}(x)\|_{L^{2}(I)}.
$$
Hence,
\begin{equation}
\label{tab}  
\|\tau_{a, b}v_{1}-\tau_{a, b}v_{2}\|_{H^{4}(I)}\leq 2\sqrt{\pi}{\cal N}_{a, \ b}l
\|v_{1}-v_{2}\|_{H^{4}(I)}.
\end{equation}
Note that the constant in the right side of (\ref{tab}) is less than one as assumed.
The Fixed Point Theorem yields the existence and
uniqueness of a function $v_{a, b}\in H^{4}(I)$, such that
$\tau_{a, b}v_{a, b}=v_{a, b}$. This is the only solution of problem (\ref{id1})
in $H^{4}(I)$ in case I) of the theorem. Suppose $v_{a, b}(x)\equiv 0$ in $I$. This will give a contradiction to our assumption that
$G_{n}F(0,x)_{n}\neq 0$ for some
$n\in {\mathbb Z}$. Let us note that in the situation of our theorem when $a>0$ the
reasoning does not use any orthogonality relations.           \hfill\lanbox

\bigskip

We proceed to demonstrating the validity of the final main proposition of the work.

\bigskip

\noindent
{\it Proof of Theorem 4.} Obviously, the limiting kernel $G(x)$ is
also a periodic function on the interval $I$ (see the reasoning of Lemma A4
of the Appendix below). Each equation (\ref{id1mi}) admits a unique solution
$u_{m}(x), \ m\in {\mathbb N}$ belonging to $H^{4}(I)$ in the situation when
$a>0$ and to $H^{4}_{0}(I)$ in the case when $a$ is trivial according to Theorem 3.
Limiting problem (\ref{id1}) has a unique solution $u(x)$, which belongs to
$H^{4}(I)$ in the case if $a>0$
and to $H^{4}_{0}(I)$ in the situation when $a=0$ due to Lemma A4 along
with Theorem 3.

We apply Fourier transform (\ref{fti}) to both sides of equations
(\ref{id1}) and (\ref{id1mi}). This gives us
\begin{equation}
\label{umnGmn}
u_{n}=\sqrt{2\pi}\frac{G_{n}\varphi_{n}}{n^{4}-a-ibn}, \quad  
u_{m, n}=\sqrt{2\pi}\frac{G_{m, n}\varphi_{m, n}}{n^{4}-a-ibn}, \quad
n\in {\mathbb Z}, \quad m\in {\mathbb N}.
\end{equation}
Here $\varphi_{n}$ and $\varphi_{m, n}$ stand for the Fourier images of $F(u(x), x)$
and $F(u_{m}(x), x)$ respectively under transform (\ref{fti}). Let us derive the
estimate from above
$$
|u_{m, n}-u_{n}|\leq \sqrt{2\pi}\Bigg\|\frac{G_{m, n}}{n^{4}-a-ibn}-
\frac{G_{n}}{n^{4}-a-ibn}\Bigg\|_{l^{\infty}}|\varphi_{n}|+
\sqrt{2\pi}\Bigg\|\frac{G_{m, n}}{n^{4}-a-ibn}\Bigg\|_{l^{\infty}}
|\varphi_{m, n}-\varphi_{n}|.
$$
Thus,
$$
\|u_{m}-u\|_{L^{2}(I)}\leq \sqrt{2\pi}\Bigg\|\frac{G_{m, n}}{n^{4}-a-ibn}-
\frac{G_{n}}{n^{4}-a-ibn}\Bigg\|_{l^{\infty}}\|F(u(x), x)\|_{L^{2}(I)}+
$$
$$
\sqrt{2\pi}\Bigg\|\frac{G_{m, n}}{n^{4}-a-ibn}\Bigg\|_{l^{\infty}}
\|F(u_{m}(x), x)-F(u(x), x)\|_{L^{2}(I)}.
$$
Let us use inequality (\ref{lk1}) of Assumption 1 to obtain
\begin{equation}
\label{Fumui}
\|F(u_{m}(x), x)-F(u(x), x)\|_{L^{2}(I)}\leq l\|u_{m}(x)-u(x)\|_{L^{2}(I)}.
\end{equation}
Note that $u_{m}(x), u(x)\in H^{4}(I)\subset L^{\infty}(I)$ due to the Sobolev
embedding. Evidently,
$$
\|u_{m}-u\|_{L^{2}(I)}\Bigg\{1-\sqrt{2\pi}l\Bigg\|\frac{G_{m, n}}{n^{4}-a-ibn}
\Bigg\|_{l^{\infty}}\Bigg\}\leq
$$
$$
\sqrt{2\pi}\Bigg\|\frac{G_{m, n}}{n^{4}-a-ibn}-
\frac{G_{n}}{n^{4}-a-ibn}\Bigg\|_{l^{\infty}}\|F(u(x), x)\|_{L^{2}(I)}.
$$
Recall bounds (\ref{cnabm}) in the situation when $a>0$ and (\ref{cn0bm})
in the case when $a=0$  to arrive at
$$
\|u_{m}-u\|_{L^{2}(I)}\leq \frac{\sqrt{2\pi}}{\varepsilon}
\Bigg\|\frac{G_{m, n}}{n^{4}-a-ibn}-\frac{G_{n}}{n^{4}-a-ibn}\Bigg\|_{l^{\infty}}
\|F(u(x), x)\|_{L^{2}(I)}.
$$
Obviously, $F(u(x), x)\in L^{2}(I)$ for $u(x)\in H^{4}(I)$ via condition
(\ref{ub1}) of Assumption 1. According to the result of Lemma A4 of the Appendix, we have
\begin{equation}
\label{umuic}
u_{m}(x)\to u(x), \quad m\to \infty
\end{equation}
in $L^{2}(I)$. Clearly,
$$
|n^{4}u_{m, n}-n^{4}u_{n}|\leq \sqrt{2\pi}\Bigg\|\frac{n^{4}G_{m, n}}{n^{4}-a-ibn}-
\frac{n^{4}G_{n}}{n^{4}-a-ibn}\Bigg\|_{l^{\infty}}|\varphi_{n}|+
$$
$$
\sqrt{2\pi}\Bigg\|\frac{n^{4}G_{m, n}}{n^{4}-a-ibn}\Bigg\|_{l^{\infty}}
|\varphi_{m, n}-\varphi_{n}|.
$$
By virtue of (\ref{Fumui}), we derive
$$
\Bigg\|\frac{d^{4}u_{m}}{dx^{4}}-\frac{d^{4}u}{dx^{4}}\Bigg\|_{L^{2}(I)}\leq
\sqrt{2\pi}\Bigg\|\frac{n^{4}G_{m, n}}{n^{4}-a-ibn}-
\frac{n^{4}G_{n}}{n^{4}-a-ibn}\Bigg\|_{l^{\infty}}\|F(u(x), x)\|_{L^{2}(I)}+
$$
$$
\sqrt{2\pi}\Bigg\|\frac{n^{4}G_{m, n}}{n^{4}-a-ibn}\Bigg\|_{l^{\infty}}l
\|u_{m}(x)-u(x)\|_{L^{2}(I)}.
$$
Let us use the result of Lemma A4 along with statement (\ref{umuic}). Hence,
$\displaystyle{\frac{d^{4}u_{m}}{dx^{4}}\to \frac{d^{4}u}{dx^{4}}}$ as
$m\to \infty$ in $L^{2}(I)$. This means that $u_{m}(x)\to u(x)$ in the $H^{4}(I)$
norm as $m\to \infty$.

Suppose that $u_{m}(x)$ is trivial in the interval $I$ for a certain
$m\in {\mathbb N}$. This will imply a contradiction to the assumption that
$G_{m, n}F(0, x)_{n}\neq 0$ for some $n\in {\mathbb Z}$. The similar
argument works for the solution $u(x)$ of limiting problem (\ref{id1}).
\hfill\lanbox

\bigskip

%%%%%%%%%%%%%%%%%%%%%%%%%%%%%%%%%%%%%%%%%%%%%%%%%%%%

\setcounter{equation}{0}

\section{Appendix}

\bigskip

Let $G(x)$ be a function, $G(x): {\mathbb R}\to {\mathbb R}$.
We designate its standard Fourier transform using the hat symbol as
\begin{equation}
\label{ft}  
\widehat{G}(p):={1\over \sqrt{2\pi}}\int_{-\infty}^{\infty}G(x)
e^{-ipx}dx, \quad p\in {\mathbb R}.
\end{equation}
Clearly,
\begin{equation}
\label{inf1}
\|\widehat{G}(p)\|_{L^{\infty}({\mathbb R})}\leq {1\over \sqrt{2\pi}}
\|G(x)\|_{L^{1}({\mathbb R})}
\end{equation}
and
$\displaystyle{G(x)={1\over \sqrt{2\pi}}\int_{-\infty}^{\infty}
\widehat{G}(q)e^{iqx}dq, \ x\in {\mathbb R}.}$ 
For the technical purposes we define the auxiliary term
\begin{equation}
\label{Na}
N_{a, \ b}:=\hbox{max}\Big\{
\Big\|{\widehat{G}(p)\over p^{4}-a-ibp}\Big\|_{L^{\infty}({\mathbb R})},
\quad
\Big\|{p^{4}\widehat{G}(p)\over p^{4}-a-ibp}\Big\|_{L^{\infty}({\mathbb R})}
\Big\}.
\end{equation}
Here $a\geq 0, \ b\in {\mathbb R}, \ b\neq 0$ are the constants.

\bigskip

\noindent
{\bf Lemma A1.} {\it Let $G(x): {\mathbb R}\to {\mathbb R}, \ G(x)\in L^{1}
({\mathbb R})$.

\medskip
  
\noindent  
a) If $a>0, \ b\in {\mathbb R}, \ b\neq 0$ then
$N_{a, \ b}<\infty$.

\medskip

\noindent
b) If $a=0, \ b\in {\mathbb R}, \ b\neq 0$ and in addition
$xG(x)\in L^{1}({\mathbb R})$ then $N_{0, \ b}<\infty$ if and only if the
orthogonality relation
\begin{equation}
\label{or1}
(G(x),1)_{L^{2}({\mathbb R})}=0 
\end{equation}
holds.}

\medskip

\noindent
{\it Proof.} First of all, it can be easily verified that in both cases a) and
b) of our lemma the boundedness of    
$\displaystyle{\frac{\widehat{G}(p)}{p^{4}-a-ibp}}$ yields the boundedness of
$\displaystyle{\frac{p^{4}\widehat{G}(p)}{p^{4}-a-ibp}}$. Evidently,
$\displaystyle{\frac{p^{4}\widehat{G}(p)}{p^{4}-a-ibp}}$ can be written as
\begin{equation}
\label{p4ghp}
\widehat{G}(p)+a\frac{\widehat{G}(p)}{p^{4}-a-ibp}+
\frac{ibp\widehat{G}(p)}{p^{4}-a-ibp}.
\end{equation}
The first term in (\ref{p4ghp}) belongs to $L^{\infty}({\mathbb R})$ due
to (\ref{inf1}) because $G(x)\in L^{1}({\mathbb R})$ as assumed.
For the third term in (\ref{p4ghp}) we have the estimate in the absolute
value as
$$
\frac{|b||p||\widehat{G}(p)|}{\sqrt{(p^{4}-a)^{2}+b^{2}p^{2}}}\leq
\|\widehat{G}(p)\|_{L^{\infty}({\mathbb R})}\leq \frac{1}{\sqrt{2\pi}}
\|G(x)\|_{L^{1}({\mathbb R})}<\infty
$$
by means of inequality (\ref{inf1}).

Hence,
$\displaystyle{\frac{\widehat{G}(p)}{p^{4}-a-ibp}\in L^{\infty}({\mathbb R})}$
implies that
$\displaystyle{\frac{p^{4}\widehat{G}(p)}{p^{4}-a-ibp}\in L^{\infty}({\mathbb R})}$
as well.

To demonstrate the validity of the result of the part a) of the lemma,
we need to consider the expression
\begin{equation}
\label{p4ghpfr}
\frac{|\widehat{G}(p)|}{\sqrt{(p^{4}-a)^{2}+b^{2}p^{2}}}.
\end{equation}  
The numerator of (\ref{p4ghpfr}) can be trivially estimated from above
using (\ref{inf1}). The denominator of (\ref{p4ghpfr}) can be bounded from
below by a finite, positive constant. Thus,
$$
\Bigg|\frac{\widehat{G}(p)}{p^{4}-a-ibp}\Bigg|\leq C\|G(x)\|_
{L^{1}({\mathbb R})}<\infty.
$$
Here and further down $C$ will denote a finite, positive constant.
Therefore, under the stated assumptions, if $a>0$ we have
$N_{a, \ b}<\infty$.

In the case of $a=0$, we write
$$
\widehat{G}(p)=\widehat{G}(0)+\int_{0}^{p}\frac{d\widehat{G}(s)}{ds}ds,
$$
such that
\begin{equation}
\label{Gpib}  
\frac{\widehat{G}(p)}{p^{4}-ibp}=\frac{\widehat{G}(0)}{p(p^{3}-ib)}+
\frac{\int_{0}^{p}\frac{d\widehat{G}(s)}{ds}ds}{p(p^{3}-ib)}.
\end{equation}
By virtue of definition (\ref{ft}) of the standard Fourier transform, we
easily obtain
$$
\Bigg|\frac{d\widehat{G}(p)}{dp}\Bigg|\leq \frac{1}{\sqrt{2\pi}}
\|xG(x)\|_{L^{1}({\mathbb R})}.
$$
Then
$$
\Bigg|\frac{\int_{0}^{p}\frac{d\widehat{G}(s)}{ds}ds}{p(p^{3}-ib)}\Bigg|\leq
\frac{\|xG(x)\|_{L^{1}({\mathbb R})}}{\sqrt{2\pi}|b|}<\infty
$$
as we assume. This means that the expression in the left side
of formula (\ref{Gpib}) is bounded if and only if $\widehat{G}(0)$
vanishes.
This is equivalent to our orthogonality condition (\ref{or1}).    \hfill\lanbox

\bigskip

Let us introduce the following technical quantities, which will help us to
study equations (\ref{id1m}).
\begin{equation}
\label{Nam}
N_{a, \ b, \ m}:=\hbox{max}\Big\{
\Big\|{\widehat{G_{m}}(p)\over p^{4}-a-ibp}\Big\|_{L^{\infty}({\mathbb R})},
\quad
\Big\|{p^{4}\widehat{G_{m}}(p)\over p^{4}-a-ibp}\Big\|_{L^{\infty}({\mathbb R})}
\Big\}
\end{equation}
with the constants $a\geq 0, \ b\in {\mathbb R}, \ b\neq 0$ and
$m\in {\mathbb N}$.
We have the following auxiliary proposition.

\bigskip

\noindent
{\bf Lemma A2.} {\it Let $m\in {\mathbb N}, \ G_{m}(x): {\mathbb R}\to
{\mathbb R}, \ G_{m}(x)\in L^{1}({\mathbb R})$ and
$G_{m}(x)\to G(x)$ in $L^{1}({\mathbb R})$ as $m\to \infty$.

\medskip

\noindent
a) If $a>0, \ b\in {\mathbb R}, \ b\neq 0$, let
\begin{equation}
\label{2rnabml}  
2\sqrt{\pi}N_{a, \ b, \ m}l\leq 1-\varepsilon
\end{equation}
for all $m\in {\mathbb N}$
with a certain fixed $0<\varepsilon<1$.

\medskip

\noindent
b) If $a=0, \ b\in {\mathbb R}, \ b\neq 0$, let
$xG_{m}(x)\in L^{1}({\mathbb R}), \ xG_{m}(x)\to xG(x)$ in $L^{1}({\mathbb R})$
as $m\to \infty$, the orthogonality condition
\begin{equation}
\label{or1m}
(G_{m}(x),1)_{L^{2}({\mathbb R})}=0, \quad m\in {\mathbb N} 
\end{equation}
is valid. Let in addition
\begin{equation}
\label{2rnabml0}  
2\sqrt{\pi}N_{0, \ b, \ m}l\leq 1-\varepsilon
\end{equation}
for all $m\in {\mathbb N}$ with some fixed $0<\varepsilon<1$.

\medskip

\noindent
Then
\begin{equation}
\label{Gmpaib}  
\frac{\widehat{G_{m}}(p)}{p^{4}-a-ibp}\to \frac{\widehat{G}(p)}{p^{4}-a-ibp},
\quad m\to \infty ,
\end{equation}
\begin{equation}
\label{Gmp2aib}  
\frac{p^{4}\widehat{G_{m}}(p)}{p^{4}-a-ibp}\to \frac{p^{4}\widehat{G}(p)}
{p^{4}-a-ibp}, \quad m\to \infty    
\end{equation}
in $L^{\infty}({\mathbb R})$, such that
\begin{equation}
\label{Gmpaibc}
\Bigg\|\frac{\widehat{G_{m}}(p)}{p^{4}-a-ibp}\Bigg\|_{L^{\infty}({\mathbb R})}\to
\Bigg\|\frac{\widehat{G}(p)}{p^{4}-a-ibp}\Bigg\|_{L^{\infty}({\mathbb R})}, \quad
m\to \infty,
\end{equation}
\begin{equation}
\label{Gmpaibcp}  
\Bigg\|\frac{p^{4}\widehat{G_{m}}(p)}{p^{4}-a-ibp}\Bigg\|_{L^{\infty}({\mathbb R})}\to
\Bigg\|\frac{p^{4}\widehat{G}(p)}{p^{4}-a-ibp}\Bigg\|_{L^{\infty}({\mathbb R})}, \quad
m\to \infty.
\end{equation}
Furthermore,
\begin{equation}
\label{Nabeps}  
2\sqrt{\pi}N_{a, \ b}l\leq 1-\varepsilon.
\end{equation}
}

\medskip

\noindent
{\it Proof.} Let us recall bound (\ref{inf1}). We easily derive
\begin{equation}
\label{GmGp}
\|\widehat{G}_{m}(p)-\widehat{G}(p)\|_{L^{\infty}({\mathbb R})}\leq
\frac{1}{\sqrt{2\pi}}\|G_{m}(x)-G(x)\|_{L^{1}({\mathbb R})}\to 0, \quad m\to \infty
\end{equation}
as assumed.
Evidently, (\ref{Gmpaibc}) and (\ref{Gmpaibcp}) 
follow from (\ref{Gmpaib}) and (\ref{Gmp2aib}) respectively
via the standard triangle inequality.

It can be trivially checked that (\ref{Gmpaib}) implies (\ref{Gmp2aib}). Let us write 
$$
\frac{p^{4}\widehat{G}_{m}(p)}{p^{4}-a-ibp}-\frac{p^{4}\widehat{G}(p)}
{p^{4}-a-ibp}=\widehat{G}_{m}(p)-\widehat{G}(p)+a\frac{\widehat{G}_{m}(p)-\widehat{G}(p)}{p^{4}-a-ibp}+
\frac{ibp[\widehat{G}_{m}(p)-\widehat{G}(p)]}{p^{4}-a-ibp},
$$
such that
$$
\Bigg\|\frac{p^{4}\widehat{G}_{m}(p)}{p^{4}-a-ibp}-\frac{p^{4}\widehat{G}(p)}{p^{4}-a-ibp}\Bigg\|_{L^{\infty}({\mathbb R})}\leq
\|\widehat{G}_{m}(p)-\widehat{G}_{m}(p)\|_{L^{\infty}({\mathbb R})}+
$$
\begin{equation}
\label{p4gmgp}
a
\Bigg\|\frac{\widehat{G}_{m}(p)}{p^{4}-a-ibp}-\frac{\widehat{G}(p)}{p^{4}-a-ibp}\Bigg\|_{L^{\infty}({\mathbb R})}+|b|
\Bigg\|\frac{p[\widehat{G}_{m}(p)-\widehat{G}(p)]}{p^{4}-a-ibp}\Bigg\|_{L^{\infty}({\mathbb R})}.
\end{equation}
The first term in the right side of  (\ref{p4gmgp}) converges to zero as $m\to \infty$ by means of (\ref{GmGp}). The second term in
the right side of  (\ref{p4gmgp}) tends to zero as $m\to \infty$ due to  (\ref{Gmpaib}). To estimate the last term in the right side of  (\ref{p4gmgp}), we observe that
$$
\frac{|p||\widehat{G}_{m}(p)-\widehat{G}(p)|}{\sqrt{(p^{4}-a)^{2}+b^{2}p^{2}}}\leq \frac{1}{|b|}
\|\widehat{G}_{m}(p)-\widehat{G}(p)\|_{L^{\infty}({\mathbb R})}.
$$
Thus,
$$
|b|\Bigg\|\frac{p[\widehat{G}_{m}(p)-\widehat{G}(p)]}{p^{4}-a-ibp}\Bigg\|_{L^{\infty}({\mathbb R})}\leq 
\|\widehat{G}_{m}(p)-\widehat{G}(p)\|_{L^{\infty}({\mathbb R})}\to 0           
$$
as $m\to \infty$ via (\ref{GmGp}). Therefore,   (\ref{Gmpaib}) yields (\ref{Gmp2aib}).

Let us demonstrate the validity of (\ref{Gmpaib}) in the case of $a>0$.
Consider the expression
\begin{equation}
\label{GmGpababs}
\frac{|\widehat{G}_{m}(p)-\widehat{G}(p)|}{\sqrt{(p^{4}-a)^{2}+b^{2}p^{2}}}. 
\end{equation}
Clearly, the denominator in  (\ref{GmGpababs}) can be bounded
from below by a positive constant and the numerator in fraction (\ref{GmGpababs}) can
be estimated from above by means of (\ref{GmGp}). Hence,
$$
\Bigg\|\frac{\widehat{G}_{m}(p)}{p^{4}-a-ibp}-\frac{\widehat{G}(p)}{p^{4}-a-ibp}
\Bigg\|_{L^{\infty}({\mathbb R})}\leq C
\|G_{m}(x)-G(x)\|_{L^{1}({\mathbb R})}\to 0
$$
as $m\to \infty$ as we assume. This means that
(\ref{Gmpaib}) holds in the situation a) of the lemma.

We turn our attention to establishing (\ref{Gmpaib}) in the case when $a$
is trivial.
In this situation orthogonality condition (\ref{or1m}) holds according to our assumption. We easily
demonstrate that the analogous relation will be valid in the limit. Obviously,
$$
|(G(x), 1)_{L^{2}({\mathbb R})}|=|(G(x)-G_{m}(x), 1)_{L^{2}({\mathbb R})}|\leq
\|G_{m}(x)-G(x)\|_{L^{1}({\mathbb R})}\to 0
$$
as $m\to \infty$ as assumed. Therefore,
\begin{equation}
\label{G1l}
(G(x), 1)_{L^{2}({\mathbb R})}=0
\end{equation}  
holds. Evidently, we can express
$$
\widehat{G}(p)=\widehat{G}(0)+\int_{0}^{p}\frac{d\widehat{G}(s)}{ds}ds, \quad
\widehat{G}_{m}(p)=\widehat{G}_{m}(0)+\int_{0}^{p}\frac{d\widehat{G}_{m}(s)}{ds}ds,
\quad m\in {\mathbb N}.
$$
Formulas (\ref{G1l}) and (\ref{or1m}) along with the definition of the standard Fourier transform (\ref{ft}) imply that
$$
\widehat{G}(0)=0, \quad \widehat{G}_{m}(0)=0, \quad m\in {\mathbb N}.
$$
We derive
\begin{equation}
\label{GmGpbub}  
\Bigg|\frac{\widehat{G}_{m}(p)}{p^{4}-ibp}-\frac{\widehat{G}(p)}{p^{4}-ibp}\Bigg|=
\Bigg|\frac{\int_{0}^{p}\Big[\frac{d\widehat{G}_{m}(s)}{ds}-\frac{d\widehat{G}
(s)}{ds}\Big]{ds}}{p^{4}-ibp}\Bigg|.
\end{equation}
The definition of the standard Fourier transform (\ref{ft}) gives us
$$
\Bigg|\frac{d\widehat{G}_{m}(p)}{dp}-\frac{d\widehat{G}(p)}{dp}\Bigg|\leq
\frac{1}{\sqrt{2\pi}}\|xG_{m}(x)-xG(x)\|_{L^{1}({\mathbb R})}.
$$
This allows us to obtain the upper bound on the right side of (\ref{GmGpbub})
equal to
$$
\frac{\|xG_{m}(x)-xG(x)\|_{L^{1}({\mathbb R})}}{\sqrt{2\pi}|b|},
$$
such that
$$
\Bigg\|\frac{\widehat{G}_{m}(p)}{p^{4}-ibp}-\frac{\widehat{G}(p)}{p^{4}-ibp}\Bigg\|_
{L^{\infty}({\mathbb R})}\leq
\frac{\|xG_{m}(x)-xG(x)\|_{L^{1}({\mathbb R})}}{\sqrt{2\pi}|b|}\to 0            
$$
as $m\to\infty$ according to the given condition. Therefore, (\ref{Gmpaib})
is valid in the situation when $a$ vanishes. Note that under the assumptions of the  lemma, we have
$$
N_{a, \ b}<\infty, \quad N_{a, \ b, \ m}<\infty, \quad m\in {\mathbb N}, \quad
a\geq 0, \quad b\in {\mathbb R}, \quad b\neq 0
$$
by means of the result of Lemma A1 above. Let us recall condition  (\ref{2rnabml})
when $a>0$ and (\ref{2rnabml0}) if $a=0$. An trivial limiting argument
using (\ref{Gmpaibc}) and (\ref{Gmpaibcp}) gives us (\ref{Nabeps}).
\hfill\lanbox 

\bigskip

Let us consider the function $G(x): I\to {\mathbb R}$, so that $G(0)=G(2\pi)$. Its
Fourier image on our finite interval is introduced  as
\begin{equation}
\label{fti}  
G_{n}:=\int_{0}^{2\pi}G(x){e^{-inx}\over \sqrt{2\pi}}dx, \quad n\in
{\mathbb Z},
\end{equation}
such that
$\displaystyle{G(x)=\sum_{n=-\infty}^{\infty}G_{n}{e^{inx}\over \sqrt{2\pi}}}$.
Obviously, the upper bound
\begin{equation}
\label{Gnub}
\|G_{n}\|_{l^{\infty}}\leq \frac{1}{\sqrt{2\pi}}\|G(x)\|_{L^{1}(I)}
\end{equation}
is valid. Evidently, if our function is continuous on the interval $I$, we
have the estimate from above
\begin{equation}
\label{cont}
\|G(x)\|_{L^{1}(I)}\leq 2\pi\|G(x)\|_{C(I)}. 
\end{equation}
Analogously to the situation on the whole real line, we define
\begin{equation}
\label{Nab}
{\cal N}_{a, \ b}:=max \Bigg\{ \Bigg\|{G_{n}\over n^{4}-a-ibn}\Bigg\|_{l^{\infty}},
\quad \Bigg\|{n^{4}G_{n}\over n^{4}-a-ibn}\Bigg\|_{l^{\infty}}\Bigg\}.
\end{equation}
Here $a\geq 0, \ b\in {\mathbb R}, \ b\neq 0$ are the constants. 
We have the following technical statement.

\bigskip

\noindent
{\bf Lemma A3.} {\it Let $\displaystyle{G(x): I\to {\mathbb R}, \
G(x)\in C(I)}$ and $G(0)=G(2\pi)$.

\medskip
  
\noindent
a) If $a>0, \ b\in {\mathbb R}, \ b\neq 0$ then ${\cal N}_{a, \ b}<\infty.$

\medskip

\noindent
b) If $a=0, \ b\in {\mathbb R}, \ b\neq 0$ then ${\cal N}_{0, \ b}<\infty$ if
and only if the orthogonality relation
\begin{equation}
\label{oc1}
(G(x),1)_{L^{2}(I)}=0
\end{equation}
holds.}

\medskip

\noindent
{\it Proof.} It can be trivially checked that in both cases a) and b) of our
lemma under the stated assumptions 
$\displaystyle{\frac{G_{n}}{n^{4}-a-ibn}\in l^{\infty}}$ implies that
$\displaystyle{\frac{n^{4}G_{n}}{n^{4}-a-ibn}\in l^{\infty}}$ as well. Clearly,
$\displaystyle{\frac{n^{4}G_{n}}{n^{4}-a-ibn}}$ can be expressed as
\begin{equation}
\label{gn4aibn}
G_{n}+a\frac{G_{n}}{n^{4}-a-ibn}+\frac{ibnG_{n}}{n^{4}-a-ibn}.
\end{equation}
The first term in (\ref{gn4aibn}) belongs to $l^{\infty}$ due to (\ref{Gnub}) and (\ref{cont})
for $G(x)\in C(I)$ as we assume. For the third term in (\ref{gn4aibn}) we have the upper bound in
the absolute value as
$$
\frac{|b||n||G_{n}|}{\sqrt{(n^{4}-a)^{2}+b^{2}n^{2}}}\leq \|G_{n}\|_{l^{\infty}}\leq \sqrt{2\pi}\|G(x)\|_{C(I)}<\infty
$$
by virtue of estimates (\ref{Gnub}) and (\ref{cont}). Therefore, the boundedness of
$\displaystyle{\frac{G_{n}}{n^{4}-a-ibn}}$ yields the boundedness of
$\displaystyle{\frac{n^{4}G_{n}}{n^{4}-a-ibn}}$.

Let us establish the statement of the part a) of our lemma. Consider
the expression
\begin{equation}
\label{Gnab}
\frac{|G_{n}|}{\sqrt{(n^{4}-a)^{2}+b^{2}n^{2}}}.  
\end{equation}
Evidently, the denominator in (\ref{Gnab}) can be easily bounded below
by a positive constant. The numerator in (\ref{Gnab}) can be trivially
estimated from above by means of (\ref{Gnub}) and (\ref{cont}).
This means that ${\cal N}_{a, \ b}<\infty$ in the case when $a>0$. It remains to demonstrate the
validity of the result of the lemma in the situation when $a=0$. Note that
\begin{equation}
\label{Gib}
\Bigg|\frac{G_{n}}{n^{4}-ibn}\Bigg|
\end{equation}
is bounded if and only if $G_{0}$ vanishes. This is equivalent to orthogonality
condition (\ref{oc1}). In such a case expression (\ref{Gib}) can be bounded
from above as
$$
\frac{|G_{n}|}{|n|\sqrt{n^{6}+b^{2}}}\leq \frac{\sqrt{2\pi}\|G(x)\|_{C(I)}}{|b|}<\infty
$$
by virtue of inequalities (\ref{Gnub}) and (\ref{cont}) under the stated assumptions.
\hfill\lanbox

\bigskip

In order to study problems (\ref{id1mi}), we introduce
\begin{equation}
\label{Nabmc}
{\cal N}_{a, \ b, \ m}:=max \Bigg\{ \Bigg\|{G_{m, n}\over n^{4}-a-ibn}\Bigg\|_
{l^{\infty}}, \quad
\Bigg\|{n^{4}G_{m, n}\over n^{4}-a-ibn}\Bigg\|_{l^{\infty}}\Bigg\},
\end{equation}
where
$a\geq 0, \ b\in {\mathbb R}, \ b\neq 0$ are the constants and $m\in {\mathbb N}$.
Let us conclude the work with the following auxiliary proposition.

\bigskip

\noindent
{\bf Lemma A4.} {\it Assume that for $m\in {\mathbb N}$ we have
$G_{m}(x): I\to {\mathbb R}, \ G_{m}(x)\in C(I), \ G_{m}(0)=G_{m}(2\pi)$
and
$G_{m}(x)\to G(x) \ in \ C(I)$ as $m\to \infty$.

\medskip

\noindent 
a) When $a>0, \ b\in {\mathbb R}, \ b\neq 0$, let
\begin{equation}
\label{cnabm}
2\sqrt{\pi}{\cal N}_{a, \ b, \ m}l\leq 1-\varepsilon
\end{equation}
for all $m\in {\mathbb N}$ with a certain fixed $0<\varepsilon<1$.

\medskip

\noindent
b) For $a=0, \ b\in {\mathbb R}, \ b\neq 0$, let the orthogonality relation
\begin{equation}
\label{or1mi}
(G_{m}(x), 1)_{L^{2}(I)}=0, \quad m\in {\mathbb N}  
\end{equation}
hold. Additionally,
\begin{equation}
\label{cn0bm}
2\sqrt{\pi}{\cal N}_{0, \ b, \ m}l\leq 1-\varepsilon
\end{equation}
for all $m\in {\mathbb N}$ with some fixed $0<\varepsilon<1$.

\medskip

\noindent
Then
\begin{equation}
\label{Gmnaib}  
\frac{G_{m, n}}{n^{4}-a-ibn}\to \frac{G_{n}}{n^{4}-a-ibn},
\quad m\to \infty ,
\end{equation}
\begin{equation}
\label{Gmn2aib}  
\frac{n^{4}G_{m, n}}{n^{4}-a-ibn}\to \frac{n^{4}G_{n}}
{n^{4}-a-ibn}, \quad m\to \infty    
\end{equation}
in $l^{\infty}$, such that
\begin{equation}
\label{Gmnaibc}
\Bigg\|\frac{G_{m, n}}{n^{4}-a-ibn}\Bigg\|_{l^{\infty}}\to
\Bigg\|\frac{G_{n}}{n^{4}-a-ibn}\Bigg\|_{l^{\infty}}, \quad
m\to \infty,
\end{equation}
\begin{equation}
\label{Gmnaibcn}  
\Bigg\|\frac{n^{4}G_{m, n}}{n^{4}-a-ibn}\Bigg\|_{l^{\infty}}\to
\Bigg\|\frac{n^{4}G_{n}}{n^{4}-a-ibn}\Bigg\|_{l^{\infty}}, \quad
m\to \infty.
\end{equation}
Moreover, the inequality
\begin{equation}
\label{cNabeps}  
2\sqrt{\pi}{\cal N}_{a, \ b}l\leq 1-\varepsilon
\end{equation}
holds.}

\medskip

{\noindent}
{\it Proof.} Note that under our conditions, the limiting kernel
function $G(x)$ will be periodic as well. Evidently,
$$
|G(0)-G(2\pi)|\leq |G(0)-G_{m}(0)|+|G_{m}(2\pi)-G(2\pi)|\leq
2\|G_{m}(x)-G(x)\|_{C(I)}\to 0, \quad m\to \infty
$$
as we assume. Thus, $G(0)=G(2\pi)$.
Using inequalities (\ref{Gnub}) and (\ref{cont}), we derive
\begin{equation}
\label{Gmnlinf}  
\|G_{m, n}-G_{n}\|_{l^{\infty}}\leq \frac{1}{\sqrt{2\pi}}\|G_{m}-G\|_{L^{1}(I)}\leq
\sqrt{2\pi}\|G_{m}-G\|_{C(I)}\to 0, \quad m\to \infty
\end{equation}
due to the given condition. It can be trivially checked
that the statements of (\ref{Gmnaib})
and (\ref{Gmn2aib}) will imply (\ref{Gmnaibc}) and (\ref{Gmnaibcn})
respectively via the triangle inequality. 

Let us demonstrate  that (\ref{Gmnaib}) yields (\ref{Gmn2aib}). We express
$$
\frac{n^{4}G_{m, n}}{n^{4}-a-ibn}-\frac{n^{4}G_{n}}{n^{4}-a-ibn}=G_{m, n}-G_{n}+a\frac{G_{m, n}-G_{n}}{n^{4}-a-ibn}+
\frac{ibn[G_{m, n}-G_{n}]}{n^{4}-a-ibn},
$$
so that
$$
\Bigg\|\frac{n^{4}G_{m, n}}{n^{4}-a-ibn}-\frac{n^{4}G_{n}}{n^{4}-a-ibn}\Bigg\|_{l^{\infty}}\leq \|G_{m, n}-G_{n}\|_{l^{\infty}}+
$$
\begin{equation}
\label{gmnn4n}
a\Bigg\|\frac{G_{m, n}-G_{n}}{n^{4}-a-ibn}\Bigg\|_{l^{\infty}}+
|b|\Bigg\|\frac{n[G_{m, n}-G_{n}]}{n^{4}-a-ibn}\Bigg\|_{l^{\infty}}.
\end{equation}
The first term in the right side of  (\ref{gmnn4n}) tends to zero as $m\to \infty$ by virtue of  (\ref{Gmnlinf}).  
The second term in the right side of  (\ref{gmnn4n})  converges to zero as $m\to \infty$ via (\ref{Gmnaib}).
Let us consider the last term in the right side of  (\ref{gmnn4n}). Clearly,
$$
\frac{|n||G_{m, n}-G_{n}|}{\sqrt{(n^{4}-a)^{2}+b^{2}n^{2}}}\leq \frac{\|G_{m, n}-G_{n}\|_{l^{\infty}}}{|b|}.
$$
Thus,
$$
|b|\Bigg\|\frac{n[G_{m, n}-G_{n}]}{n^{4}-a-ibn}\Bigg\|_{l^{\infty}}\leq \|G_{m, n}-G_{n}\|_{l^{\infty}}\to 0
$$
as $m\to \infty$ according to (\ref{Gmnlinf}). Therefore, (\ref{Gmnaib}) implies (\ref{Gmn2aib}).

First we establish the validity of (\ref{Gmnaib}) in the situation when
$a>0$. Let us consider
\begin{equation}
\label{Gmnconv}
\frac{|G_{m, n}-G_{n}|}{\sqrt{(n^{4}-a)^{2}+b^{2}n^{2}}}.
\end{equation}  
Obviously, the denominator in fraction (\ref{Gmnconv}) can be estimated from below by a
positive constant. The numerator in (\ref{Gmnconv}) is bounded from above according to
(\ref{Gmnlinf}). This yields (\ref{Gmnaib}) for $a>0$.

Finally, we establish (\ref{Gmnaib}) in the case when $a=0$. Let us recall
orthogonality condition (\ref{or1mi}).
It can be trivially checked that the analogous relation holds in the limit.
Evidently,
$$
|(G(x), 1)_{L^{2}(I)}|=|(G(x)-G_{m}(x), 1)_{L^{2}(I)}|\leq 2\pi\|G_{m}(x)-G(x)\|_
{C(I)}\to 0, \quad m\to \infty
$$
as we assume. Hence,
$$
(G(x), 1)_{L^{2}(I)}=0.
$$
This is equivalent to $G_{0}=0$. Note that $G_{m, 0}=0, \ m\in {\mathbb N}$ 
via orthogonality condition (\ref{or1mi}). Let us use
(\ref{Gmnlinf}) to obtain
$$
\Bigg|\frac{G_{m, n}-G_{n}}{n^{4}-ibn}\Bigg|\leq \frac{\sqrt{2\pi}\|G_{m}(x)-G(x)\|_
{C(I)}}{|b|}.
$$
The norm in the right side of this inequality tends to zero as
$m\to \infty$.  This means that (\ref{Gmnaib}) holds in the situation when $a$ is trivial as well.
Clearly, under the stated assumptions we have
$$
{\cal N}_{a, \ b}<\infty, \quad {\cal N}_{a, \ b, \ m}<\infty, \quad m\in
{\mathbb N}, \quad a\geq 0, \quad b\in {\mathbb R}, \quad b\neq 0
$$
by virtue of the result of Lemma A3 above. Let us recall bounds (\ref{cnabm})
when $a>0$ and (\ref{cn0bm}) for $a=0$. An easy limiting
argument based on (\ref{Gmnaibc}) and (\ref{Gmnaibcn}) gives us (\ref{cNabeps}).
\hfill\lanbox

\bigskip

%%%%%%%%%%%%%%%%%%%%%%%%%%%%%%%%%%%%%%%%%%%%%%%%%%%%%%%%%%%%%%%%

\section*{Acknowledgements}

The author is grateful to Israel Michael Sigal for the partial support
by the NSERC grant NA 7901.

\end{document}